\documentclass{article}

\usepackage{amssymb}

\usepackage{graphicx}

\newtheorem{corollary}{Corollary}

\newtheorem{proposition}{Proposition}

\begin{document}
	
\title{On the solution manifold of a differential equation with a delay which has a zero} 
	
\author{Hans-Otto Walther}
	
	
	
\maketitle
	
	
\begin{abstract}
 For a differential equation with a state-dependent delay we show that the associated solution manifold $X_f$ of codimension 1 in the space $C^1([-r,0],\mathbb{R})$ is an almost graph over a hyperplane, which implies that $X_f$ is diffeomorphic to the hyperplane. For the case considered previous results only provide a covering by 2 almost graphs.
\end{abstract}
	
\bigskip
	
	

\section{Introduction}



Let $r>0$ be given, choose a norm on $\mathbb{R}$, and let $C_n=C([-r,0],\mathbb{R}^n)$ and $C^1_n=C^1([-r,0],\mathbb{R}^n)$ denote the Banach spaces of continuous and continuously differentiable functions $[-r,0]\to\mathbb{R}^n$, respectively, with the norms given by $|\phi|_C=\max_{–r\le t\le0}|\phi(t)|$ and  $|\phi|=|\phi|_C+|\phi'|_C$. For a delay differential equation
$$
x'(t)=f(x_t)
$$
with a vector-valued functional $f:C^1_n\supset U\to\mathbb{R}^n$ and with the solution segment $x_t\in U$ defined as $x_t(s)=x(t+s)$, the associated {\it solution manifold} is the set
$$
X_f=\{\phi\in U:\phi'(0)=f(\phi)\}.
$$
Assume that $f$ is continuously differentiable and 

(e) {\it each derivative $Df(\phi):C^1_n\to\mathbb{R}^n$, $\phi\in U$, has a linear extension $D_ef(\phi):C_n\to\mathbb{R}^n$ so that the map 
$$
U\times C_n\ni(\phi,\chi)\mapsto D_ef(\phi)\chi\in\mathbb{R}^n
$$
is continuous.} 

The extension property (e) is a relative of the notion of being {\it almost Fr\'echet differentiable} from \cite{M-PNP} and can be verified for a variety of differential equations with state-dependent delay. - Under the said, mild conditions a non-empty set $X_f$ is a continuously differentiable submanifold of codimension $n$ in the Banach space $C^1_n$, and it is on this manifold that the initial value problem associated with Eq. (f) is well-posed, with solutions which are  continuously differentiable with respect to initial data \cite{W1,HKWW}.

The present note continues the description of solution manifolds initiated in  \cite{W6,W7}. In \cite{W6} we saw that in case $f$ satisfies a condition 
which in  examples with explicit delays corresponds to all of these delays being bounded away from zero the solution manifold is a graph over the closed subspace
$$
X_0=\{\phi\in C^1:\phi'(0)=0\},
$$
which also is the solution manifold for $f=0$.
An example of a scalar equation  with a single state-dependent delay which is positive but not bounded away from zero shows that in general solution manifolds do not admit any graph representation \cite[Section 3]{W6}. However, the main result from \cite{W6} guarantees that for a reasonably large class of systems with explicit discrete state-dependent delays which are all positive the solution manifolds are nearly as simple as graphs: They are
{\it almost graphs} over $X_0$, in the terminology introduced in \cite{W6,W7}. 

Let us recall from \cite{W7} that given
a continuously differentiable submanifold $X$ of a Banach space $E$ and a closed subspace $H$ with a closed complementary space, 

(i) $X$ is called a graph (over $H$) if there are a closed complementary space $Q$ for $H$ and a continuously differentiable map $\gamma:H\supset dom\to Q$ with
$$
X=\{\zeta+\gamma(\zeta)\in E:\zeta\in dom\},
$$

and that

(ii) $X$ is called an {\it almost graph} (over $H$)
if there is a continuously differentiable map $\alpha:H\supset dom\to E$ with
\begin{eqnarray}
\alpha(\zeta) & = & 0\quad\mbox{on}\quad dom\,\cap\,X,\nonumber\\
\alpha(\zeta) & \in & E\setminus H\quad\mbox{on}\quad dom\setminus \,X,\nonumber
\end{eqnarray}
so that the map 
$$
dom\ni\zeta\mapsto\zeta+\alpha(\zeta)\in E
$$ 
defines a diffeomorphism onto $X$.

Furthermore,

(iii) a diffeomorphism $A$ from an open neighbourhood ${\mathcal O}$ of $X$ in $E$ onto an open subset of $E$ is called an {\it almost graph diffeomorphism} (associated with $X$ and $H$) if
$$
A(X)\subset H
$$
and
$$
A(\zeta)=\zeta\quad\mbox{on}\quad X\cap\,H.
$$

In \cite[Section 1]{W7} it is verified that in case there is an almost graph diffeomorphism $A$ associated with $X$ and $H$ the submanifold $X$ is an almost graph over $H$. 

An example of an almost graph in finite dimension is the unit circle in the plane without the point on top. The inverse of the stereographic projection onto the real line serves as the map $\zeta\mapsto\zeta+\alpha(\zeta)$ in Part (ii) of the previous definition. 

We return to the results in \cite{W6} and \cite{W7}. It is not difficult to see that the approach used in \cite{W6} for an almost graph representation of a solution manifold fails in case some of the  delays in the system considered have zeros. For the solution manifold of such a system, with $k$ delays some of which have zeros, the approach can be used, however, in order to obtain a finite atlas of manifold charts whose domains are almost graphs over $X_0$. This is achieved in \cite{W7}, with the size of the atlas independent of the number of equations in the system considered, solely determined by the zerosets of the delays, and not exceeding $2^k$.

The immediate question with regard to the topological properties of these solution manifolds is whether the size of the atlas found in \cite{W7} is minimal. The result of the present note shows that for a prototype of the systems studied in \cite{W6,W7} this is not the case: The entire solution manifolds of the prototype equation is in fact an almost graph over $X_0$. The proof relies on a major modification of the approach taken in \cite{W7}.    

The prototype equation belongs to the simplest cases of the systems studied in \cite{W6,W7} which are scalar equations with a single delay ($k$=1) and have the form
\begin{equation}
x'(t)=g(x(t-d(Lx_t)))
\end{equation}
with a continuously differentiable map $g:\mathbb{R}\to\mathbb{R}$, a continuously differentiable {\it delay function} $d:F\to[0,r]$ defined on a finite-dimensional topological vectorspace $F$, and a surjective continuous linear map $L:C\to F$. We abbreviate $C=C_1$ and $C^1=C^1_1$. For $f:C^1\to\mathbb{R}$ given by
$$
f(\phi=g(\phi(-d(L\phi)))
$$ 
Eq. (1) takes the form $x'(t)=f(x_t)$. Proposition 2.1 from \cite{W7} applies and shows that $f$ is continuously differentiable with property (e). Proposition 2.3 from \cite{W7} yields that the associated solution manifold 
$$
X_f=\{\phi\in C^1:\phi'(0)=g(\phi(-d(L\phi)))\}.
$$
is non-empty, and it follows that it is a continuously differentiable submanifold of codimension $1$ in $C^1$.

If $d(\xi)>0$ everywhere than we know from \cite{W6} that $X_f$ is an almost graph over $X_0=\{\chi\in C^1:\chi'(0)=0\}$. If $d$ has zeros then the result of \cite{W7} yields an atlas of $2=2^1$ manifold charts whose domains are almost graphs over $X_0$.

In \cite{W8} an adaptation of the approach from \cite{W6,W7} to Eq. (1) with a linear map $L$ for which, losely spoken,  $L\phi$ does not depend on $\phi(0)$ is used to prove that the associated solution manifold is an almost graph over $X_0$, no matter whether $d$ has zeros or not.

In the sequel we consider a prototype for the remaining, critical cases, namely, Eq. (1) for
$$
F=\mathbb{R},\quad L\phi=\phi(0)\quad\mbox{for all}\quad\phi\in C,
$$
and for
$$
d\quad\mbox{with a single zero}\quad\eta_0\in\mathbb{R}.
$$
We find an almost graph diffeomorphism $A:C^1\to C^1$ which maps $X_f$ onto $X_0$.

A part of the construction of the diffeomorphism $A$ uses the technique developed in \cite{W6,W7}. For another part it was helpful to have in mind an idea of Krisztin \cite{K} which yields graph representations of solution manifolds from bounds on extended derivatives as in property (e), compare the proof of Lemma 1 in \cite{KR}.

The assumption that $d$ has a single zero is for simplicity and may be relaxed in future work.

\section{Preliminaries} 

Differentiable maps are always defined on open subsets of Banach spaces or Banach manifolds.

Differentiation $\partial:C^1\to C$, $\partial\phi=\phi'$, is linear and continuous, and the evaluation map $ev:C\times[-r,0]\ni(\phi,t)\mapsto\phi(t)\in\mathbb{R}$ is continuous but not locally Lipschitz continuous. The composition 
$$
ev(\cdot,0)\circ\partial:C^1\ni\phi\mapsto ev(\partial\phi,0)\in\mathbb{R}
$$
is linear and continuous.

The restriction $ev_1$ of $ev$ to $C^1\times(-r,0)$ is continuously differentiable with 
$$
D\,ev_1(\phi,t)(\chi,s)=\chi(t)+\phi'(t)\,s,
$$
and the composition 
$$
h:C^1\ni\phi\mapsto ev(\phi,-d(L\phi))\in\mathbb{R}
$$
is continuously differentiable with
$$
Dh(\phi)\chi=\chi(-d(L\phi))-\phi'(-d(L\phi))d'(L\phi)L\chi,
$$ 
see Part  2.1 of the proof of \cite[Proposition 2.1]{W7}. 

In Section 1 we quoted \cite[Proposition 2.3]{W7} for $X_f\neq\emptyset$. Using the choice of $L$ this also follows directly from the fact that for every $\xi\in\mathbb{R}$ each $\phi\in C^1$ with $L(\phi)=\phi(0)=\xi$, $\phi(-d(\xi))=\xi$, and $\phi'(0)=g(\xi)$ belongs to $X_f$.

The tangent space $T_{\phi}X_f$ of the solution manifold $X_f$ at $\phi\in X_f$ consists of the vectors $c'(0)=Dc(0)1$ of differentiable curves $c:I\to C^1$ with $0\in I\subset\mathbb{R}$, $c(0)=\phi$, $c(I)\subset X_f$. We have
\begin{eqnarray}
T_{\phi}X_f & = & \{\chi\in C^1:\chi'(0)=Df(\phi)\chi\}\nonumber\\
& = & \{\chi\in C^1:\nonumber\\
& & \chi'(0)=g'(\phi(-d(L\phi)))[\chi(-d(L\phi))-\phi'(-d(L\phi))d'(L\phi)L\chi]\}.\nonumber
\end{eqnarray}

\section{A map $A$ taking $X_f$ into $X_0$}

As $d$ is minimal at $\eta_0$ we have $d'(\eta_0)=0$.

Notice that for $\phi\in C$ with $L\phi=\eta_0$,
\begin{equation}
\phi(-d(L\phi))=\phi(-d(\eta_0))=\phi(0)=L\phi=\eta_0.
\end{equation}

For reals $\eta$ we introduce the continuous linear maps
$$
L_{\eta}:C\to\mathbb{R}
$$
given by $L_{\eta}\phi=\phi(-d(\eta))$. 

In order to develop a bit of intuition about the shape of $X_f$ observe that the sets
\begin{eqnarray}
X_{f\eta} & = & X_f\cap L^{-1}(\eta)\nonumber\\
& = & \{\phi\in C^1:\,L\phi=\eta\,\,\mbox{and}\,\,\phi'(0)=g(\phi(-d(L\phi)))
\}\nonumber
\\
& = & \{\phi\in C^1:\,L\phi=\eta\,\,\mbox{and}\,\,\phi'(0)=g(\phi(-d(\eta)))
\}\nonumber
\end{eqnarray}
are mutually disjoint and decompose $X_f$, and that
$$
X_{f\eta_0}=\{\phi\in C^1:\,L\phi=\eta_0\,\,\mbox{and}\,\,\phi'(0)=g(\eta_0)\}
$$
(with Eq. (2) for $L\phi=\eta_0$) is a closed affine subspace of codimension 2 in $C^1$.


Choose $\rho>0$ and set
\begin{eqnarray}
c & = & \max_{|\xi-\eta_0|\le\rho}|g(\xi)|+\max_{|\xi-\eta_0|\le\rho}|g'(\xi)|,\nonumber\\
c_{\ast} & = & \frac{\rho}{4(c+1)(\rho+3)}.\nonumber
\end{eqnarray}

The map which we are going to construct relies on vectors $\psi_{\eta}\in C^1$ which are transversal to the solution manifold at points in $X_{f\eta}$. We begin with the case $\eta=\eta_0$ and choose $\psi_{\eta_0}\in C^1$ with the properties
$$
\psi_{\eta_0}(0)=0,\,\psi_{\eta_0}'(0)=1,\,|\psi_{\eta_0}|_C\le c_{\ast}.
$$
For each $\phi\in X_f$ with $L\phi=\eta_0$ (or, for each $\phi\in X_{f\eta_0}$) we have
\begin{eqnarray}
\psi_{\eta_0}'(0) & = & 1>c\,c_{\ast}\ge|g(\eta_0)||\psi_{\eta_0}|_C\nonumber\\
& = & |g(\phi(-d(L\phi)))||\psi_{\eta_0}|_C\ge  |g(\phi(-d(L\phi)))\psi_{\eta_0}(-d(L\phi))|\nonumber\\
& = & |g(\phi(-d(L\phi)))\psi_{\eta_0}(-d(L\phi))-\phi'(-d(L\phi))d'(L\phi)L\psi_{\eta_0}]|\nonumber\\
& & \mbox{(with}\quad d'(L\phi)=d'(\eta_0)=0)\nonumber
\end{eqnarray}
which means
$$
\psi_{\eta_0}\in C^1\setminus T_{\phi}X_f.
$$

\begin{proposition}
There exists a continuously differentiable map
$$
\mathbb{R}\setminus\{\eta_0\}\ni\eta\mapsto\psi_{\eta}\in C^1
$$	
so that for every $\eta\in\mathbb{R}\setminus\{\eta_0\}$ we have
$$
\psi_{\eta}(t)=0\,\,\mbox{on}\,\,[-r,-d(\eta)]\cup\{0\},\,\psi_{\eta}'(0)=1,\,\,
|\psi_{\eta}|_C\le c_{\ast},
$$
and for all $\phi\in X_f$ with $L\phi=\eta$,
$$
\psi_{\eta}\in C^1\setminus T_{\phi}X_f.
$$
\end{proposition}

{\bf Proof.} 1. For each $z\in[-r,0)$ choose $\psi_z\in C^1$ with $\psi_z(t)=0$ on $[-r,z]\cup\{0\}$, $\psi'(0)=1$, and $|\psi_z|_C\le c_{\ast}$. Then proceed as in the proof of 
\cite[Proposition 4.1]{W7}, with ${\mathcal F}=\mathbb{R},{\mathcal W}=\mathbb{R}\setminus\{\eta_0\},\lambda=L$,  and construct the desired map $\mathbb{R}\setminus\{\eta_0\}\to C^1$ from a sequence of maps $\psi_{z_m}$, $m\in\mathbb{N}$, with $z_m\to\min d=0$ as $m\to\infty$. Observe that $|\psi_{\eta}|_C\le c_{\ast}$ is achieved.

2. For $\eta\in\mathbb{R}\setminus\{\eta_0\}$ and $\phi\in X_f$ with $L\phi=\eta$ the function $\psi_{\eta}$
does not satisfy the equation characterizing the tangent space $T_{\phi}X_f$, due to $\psi_{\eta}'(0)=1$ and $\psi_{\eta}(-d(L\phi))=\psi_{\eta}(-d(\eta))=0=\psi_{\eta}(0)=L\psi_{\eta}$. $\Box$

The map from Proposition 1 has no continuous extension to $\mathbb{R}$. 
Nevertheless, for all $\eta\in\mathbb{R}$,
\begin{equation}
L\psi_{\eta}=\psi_{\eta}(0)=0.
\end{equation}
Also, for each $\eta\in\mathbb{R}$,
\begin{equation}
L_{\eta}\psi_{\eta}=0,
\end{equation}
which in case $\eta=\eta_0$ holds with $L_{\eta_0}=L$, and,
\begin{equation}
\psi_{\eta}'(0)=1. 
\end{equation}
It follows from Eq. (5) that for all $\eta\in\mathbb{R}$,
$$
C^1=X_0\oplus\mathbb{R}\psi_{\eta},
$$
and the continuous linear projection $P_{\eta}:C^1\to C^1$
along $\mathbb{R}\psi_{\eta}$ onto $X_0$  is given by
$$
P_{\eta}\phi=\phi-\phi'(0)\psi_{\eta}.
$$

Now we are ready for the definition of the map $A:C^1\to C^1$ which in the next section will be shown to be an almost graph diffeomorphism associated with  $X_f$ and $X_0$. Let 
$$
a:\mathbb{R}\to[0,1]
$$
be a continuously differentiable map  with $a(\xi)=1$ for $|\xi-\eta_0|\le \frac{\rho}{2}$, $a(\xi)=0$ for $|\xi-\eta_0|\ge\rho$, and $|a'(\xi)|\le\frac{3}{\rho}$ for $\frac{\rho}{2}\le|\xi-\eta_0|\le\rho$.

The maps $A_{\rho/2}:C^1_{\rho/2}\to C^1$ and $A_+:C^1_+\to C^1$
given by
\begin{eqnarray}
C^1_{\rho/2} & = & \{\phi\in C^1:|\phi(-d(L\phi))-\eta_0|<\rho/2\},\nonumber\\
A_{\rho/2}(\phi) & = & \phi-g(\tau)\psi_{\eta_0},\nonumber\\
C^1_+ & = & \{\phi\in C^1:d(L\phi)>0\},\nonumber\\
A_+(\phi) & = & \phi-g(\tau)[a(\tau)\psi_{\eta_0}+(1-a(\tau))\psi_{\eta}],\nonumber
\end{eqnarray}
with
$$ \tau=L_\eta\phi=ev(\phi,-d(L\phi)),\quad\eta=L\phi,
$$
are continuously differentiable. On the intersection of their domains they coincide: For $\phi\in C^1_{\rho/2}\cap C^1_+$ we have  
$$
|\tau-\eta_0|=|\phi(-d(L\phi))-\eta_0|<\frac{\rho}{2},
$$
which yields $a(\tau)=1$, and thereby, $A_{\rho/2}(\phi)=A_+(\phi)$. 

Also,    
$$
C^1=C^1_{\rho/2}\cup C^1_+,
$$ 
since for $\phi\in C^1\setminus C^1_+$,  $d(L\phi)=0$, hence  $L\phi=\eta_0$, and due to Eq. (2), $\phi(-d(L\phi))=\eta_0$, which means
$|\phi(-d(L\phi))-\eta_0|=0$, or $\phi\in C^1_{\rho/2}$.

It follows that $A_{\rho/2}$ and $A_+$ define a continuously differentiable map $A:C^1\to C^1$.

Using Eq. (5) we infer that for every $\phi\in X_f$,
$$
(A(\phi))'(0)=\phi'(0)-g(\tau)\cdot1=\phi'(0)-g(\phi(-d(L\phi)))=0,
$$
or, $A(\phi)\in X_0$. Hence
$$
A(X_f)\subset X_0.
$$
Notice also that due to Eq. (3),
\begin{equation}
LA(\phi)=L\phi\quad\mbox{for all}\quad\phi\in C^1.
\end{equation}

\section{$A$ is an almost graph diffeomorphism}

In order to find the inverse of $A$ we first consider $\phi\in C^1_+\cap X_f$ and $\chi=A(\phi)=A_+(\phi)\in X_0$, and compare
$\tau=L_{\eta}\phi$ where $\eta=L\phi$ to $\sigma=L_{\hat{\eta}}\chi$
where $\hat{\eta}=L\chi$. By Eq. (6), $\hat{\eta}=\eta$, hence
\begin{eqnarray}
\sigma & = & L_{\eta}\chi=L_{\eta}A_+(\phi)\nonumber\\
& = & L_{\eta}\phi-g(\tau)[a(\tau)L_{\eta}\phi_{\eta_0}+(1-a(\tau))L_{\eta}\phi_{\eta}]\nonumber\\
& = & \tau-g(\tau)a(\tau)L_{\eta}\phi_{\eta_0}\quad\mbox{(with Eq. (4))}.
\end{eqnarray}
For every $\eta\in\mathbb{R}$ the map
$$
h_{\eta}:\mathbb{R}\ni\tau\mapsto \tau-(ga)(\tau)L_{\eta}\psi_{\eta_0}\in\mathbb{R}
$$
is continuously differentiable. 

\begin{proposition}
Every map $h_{\eta}$, $\eta\in\mathbb{R}$, is a bijection, the map 
$T:\mathbb{R}^2\to\mathbb{R}$ given by $T(\eta,\sigma)=(h_{\eta})^{-1}(\sigma)$ is
continuously differentiable, and for all  $(\eta,\sigma,\tau)\in\mathbb{R}^3$ Eq. (7) is equivalent to $\tau=T(\eta,\sigma)$.
\end{proposition}

{\bf Proof.} 1. For every $\eta\in\mathbb{R}$ the map
$h_{\eta}$ satisfies $h_{\eta}(\tau)=\tau$ for $|\tau-\eta_0|\ge\rho$, and for
$|\tau-\eta_0|\le\rho$,
\begin{eqnarray}
h_{\eta}'(\tau)& = & 1-(g'(\tau)a(\tau)+g(\tau)a'(\tau))L_{\eta}\psi_{\eta_0}\nonumber\\
& \ge & 1-(c+c(3/\rho))c_{\ast}>0.\nonumber
\end{eqnarray}
It follows that $h_{\eta}$ is bijective. 

2. The map  
$$
F:\mathbb{R}^3\ni(\eta,\sigma,\tau)\mapsto\sigma-h_{\eta}(\tau)\in\mathbb{R}
$$
is continuously differentiable (with $L_{\eta}\psi_{\eta_0}=\psi_{\eta_0}(-d(\eta))$). For each
$(\eta,\sigma,\tau)\in\mathbb{R}^3$ Eq. (7) and the relations $F(\eta,\sigma,\tau)=0$ and $\tau=(h_{\eta})^{-1}(\sigma)=T(\eta,\sigma)$
are equivalent. 
For $|\tau-\eta_0|\le\rho$ we have
$$
\partial_3F(\eta,\sigma,\tau)=-1+(ga)'(\tau)L_{\eta}\psi_{\eta_0}\le-1+(c+c(3/\rho))c_{\ast}<0,
$$
and
$$
\partial_3F(\eta,\sigma,\tau)=-1\neq0
$$
for $|\tau-\eta_0|\ge\rho$. Applications of the Implicit Function Theorem to the zeroset of $F$ show that the map $T$ is locally given by continuously differentiable maps. $\Box$

\medskip

For the open subsets $C^1_{\rho/4}=\{\chi\in C^1:|\chi(-d(L\chi))-\eta_0|<\rho/4\}$ and 
$C^1_+$ of the space $C^1$ we have
$$
C^1=C^1_{\rho/4}\cup C^1_+
$$
since for $\chi\in C^1\setminus C^1_+$,  $d(L\chi)=0$, hence  $L\chi=\eta_0$, and due to Eq. (2) $\chi(-d(L\chi))=\eta_0$, which means
$|\chi(-d(L\chi))-\eta_0|=0$. The maps $B_{\rho/4}:C^1_{\rho/4}\to C^1$ and $B_+:C^1_+\to C^1$ given by
\begin{eqnarray}
B_{\rho/4}(\chi) & = & \chi+g(\tau)\psi_{\eta_0},\nonumber\\
B_+(\chi) & = & \chi+g(\tau)[a(\tau)\psi_{\eta_0}+(1-a(\tau))\psi_{\eta}],\nonumber
\end{eqnarray}
with
$$ \tau=T(\eta,\sigma),\quad\eta=L\chi,\quad\sigma=L_\eta\chi=ev(\chi,-d(L\chi))
$$
are continuously differentiable. On the intersection $C^1_{\rho/4}\cap C^1_+$ both maps $B_{\rho/4}$ and $B_+$ coincide. In order to verify this we need to know $a(\tau)=1$ for 
$$
\tau=T(\eta,\sigma),\quad\eta=L\chi,\quad\chi\in C^1,\quad \sigma=L_{\eta}\chi
$$ 
with 
$|\chi(-d(L\chi))-\eta_0|<\rho/4$ and $d(L\chi)>0$.  The equation $a(\tau)=1$
holds provided $|\tau-\eta_0|\le\rho/2$, which follows from
$$
|\tau-\eta_0|=|\sigma+(ga)(\tau)L_{\eta}\psi_{\eta_0}-\eta_0|\le
\frac{\rho}{4}+\max_{\xi\in\mathbb{R}}|(ga)(\xi)|\cdot c_{\ast}
$$
in combination with
\begin{eqnarray}
(ga)(\tau) & = & 0\quad\mbox{for}\quad|\tau-\eta_0|\ge\rho,\nonumber\\
|(ga)(\tau)| & \le & c\quad\mbox{for}\quad|\tau-\eta_0|\le\rho,\nonumber\\
c\,c_{\ast} & < & \rho/4.\nonumber
\end{eqnarray}

So $B_{\rho/4}$ and $B_+$ define a
continuously differentiable map $B:C^1\to C^1$.

Observe that due to Eq. (3),
\begin{equation}
LB(\chi)=L\chi\quad\mbox{for all}\quad\chi\in C^1.
\end{equation}

We have
$$
B(X_0)\subset X_f.
$$
Proof of this: For $\chi\in X_0$ let $\phi=B(\chi)$ and $\eta=L\chi$, $\sigma=L_{\eta}\chi$, and $\tau=T(\eta,\sigma)=h_{\eta}^{-1}(\sigma)$.
Then 
$$
\sigma=h_{\eta}(\tau)=\tau-(ga)(\tau)L_{\eta}\psi_{\eta_0}.
$$
Using $\chi'(0)=0$, Eq. (5), and the preceding equation we get
$$
\phi'(0)=(B(\chi))'(0)=0+g(\tau)\cdot1=g(\tau)
$$
with
\begin{eqnarray}
\tau & = & \sigma+(ga)(\tau)L_{\eta}\psi_{\eta_0}\nonumber\\
& = & L_{\eta}\chi+g(\tau)[a(\tau)L_{\eta}\psi_{\eta_0}+(1-a(\tau))L_{\eta}\psi_{\eta}]\quad\mbox{(with Eq. (4))}\nonumber\\
& = & L_{\eta}B(\chi)=L_{\eta}\phi.\nonumber
\end{eqnarray}
Eq. (8) yields $\eta=L\chi=L\phi$, and we obtain
$$
\phi'(0)=g(\tau)=g(L_{\eta}\phi)=g(\phi(-d(\eta)))=g(\phi(-d(L\phi))),
$$
or, $B(\chi)=\phi\in X_f$.

\begin{proposition}
$B(A(\phi))=\phi$ for all $\phi\in C^1$.
\end{proposition}

{\bf Proof.} 1. The case $d(L\phi)>0$. Consider 
$$
\chi=A(\phi)=A_+(\phi)=\phi-g(\tau)[a(\tau)\psi_{\eta_0}+(1-a(\tau))\psi_{L\phi}]
$$ 
with
$$
\tau=L_{\eta}\phi,\quad\eta=L\phi.
$$
From Eq. (6) we infer
\begin{equation}
L\chi=L\phi,
\end{equation}
hence $d(L\chi)=d(L\phi)>0$, and thereby 
$$
B(\chi)=B_+(\chi)=\chi+g(\hat{\tau})[a(\hat{\tau})\psi_{\eta_0}+(1-a(\hat{\tau}))\psi_{L\chi}]
$$ 
with
$$
\hat{\tau}=T(\hat{\eta},\hat{\sigma})=h_{\hat{\eta}}^{-1}(\hat{\sigma}),\quad\hat{\eta}=L\chi,\quad\hat{\sigma}=L_{\hat{\eta}}\chi.
$$
It follows that
\begin{eqnarray}
B(\chi) & = & B(A(\phi))=\{ \phi-g(\tau)[a(\tau)\psi_{\eta_0}+(1-a(\tau))\psi_{L\phi}]\}\nonumber\\
& & +g(\hat{\tau})[a(\hat{\tau})\psi_{\eta_0}+(1-a(\hat{\tau}))\psi_{L\chi}].\nonumber
\end{eqnarray}
Using the preceding equation in combination with  Eq. (9) we obtain $B(A(\phi))=\phi$ provided we have $\hat{\tau}=\tau$. Proof of this:
By Eq. (9), $\hat{\eta}=L\chi=L\phi=\eta$. We get
\begin{eqnarray}
h_{\eta}(\hat{\tau}) & = & h_{\hat{\eta}}(\hat{\tau})=\hat{\sigma}=L_{\hat{\eta}}\chi=L_{\eta}\chi=L_{\eta}A(\phi)\nonumber\\
& = & L_{\eta}\phi-g(\tau)[a(\tau)L_{\eta}\psi_{\eta_0}+(1-a(\tau))L_{\eta}\psi_{\eta}]\nonumber\\
& = & L_{\eta}\phi-g(\tau)a(\tau)L_{\eta}\psi_{\eta_0}\quad\mbox{(with Eq. (4))}\nonumber\\
& = & \tau-(ga)(\tau)L_{\eta}\psi_{\eta_0}\nonumber\\
& = & h_{\eta}(\tau),\nonumber
\end{eqnarray}
and the injectivity of $h_{\eta}$ yields $\hat{\tau}=\tau$.

2. The case $d(L\phi)=0$. Then $\phi(0)=L\phi=\eta_0$. Choose a sequence
of  points $\phi_j\in C^1$, $j\in\mathbb{N}$, with $\phi_j(0)\neq\eta_0$ for all $j\in\mathbb{N}$ and $\phi_j\to\phi$ in $C^1$ as $j\to\infty$.
For all $j\in\mathbb{N}$, $B(A(\phi_j))=\phi_j$, due to Part 1 of the proof, and continuity yields $B(A(\phi))=\phi$. $\Box$

\begin{proposition}
	$A(B(\chi))=\chi$ for all $\chi\in C^1$.
\end{proposition}

{\bf Proof.} 1. The case $d(L\chi)>0$. Then $\chi\in C^1_+$ and
$$
B(\chi)=B_+(\chi)=\chi+g(\tau)[a(\tau)\psi_{\eta_0}+(1-a(\tau))\psi_{\eta}]
$$
with $\tau=T(\eta,\sigma)=h_{\eta}^{-1}(\sigma)$, $\eta=L\chi$, $\sigma=L_{\eta}\chi$. Set $\phi=B(\chi)$. By Eq. (8), $L\phi=L\chi=\eta$. Hence $d(L\phi)=d(L\chi)>0$, or $\phi\in C_+^1$, and
$$
A(\phi) = A_+(\phi)=\phi-g(\hat{\tau})[a(\hat{\tau})\psi_{\eta_0}+(1-a(\hat{\tau}))\psi_{\hat{\eta}}]
$$
with
$$
\hat{\tau}=L_{\hat{\eta}}\phi,\quad\hat{\eta}=L\phi\quad(=\eta).
$$
It follows that
\begin{eqnarray}
A(B(\chi)) & = & A(\phi)=\{\chi+g(\tau)[a(\tau)\psi_{\eta_0}+(1-a(\tau))\psi_{\eta}]\}\nonumber\\
& & -g(\hat{\tau}[a(\hat{\tau})\psi_{\eta_0}+(1-a(\hat{\tau}))\psi_{\hat{\eta}}]\nonumber\\
 & = & \{\chi+g(\tau)[a(\tau)\psi_{\eta_0}+(1-a(\tau))\psi_{\eta}]\}\nonumber\\
& & -g(\hat{\tau}[a(\hat{\tau})\psi_{\eta_0}+(1-a(\hat{\tau}))\psi_{\eta}]
\quad\mbox{(with}\quad\hat{\eta}=\eta),\nonumber
\end{eqnarray}
and for $A(B(\chi))=\chi$ it remains to show $\hat{\tau}=\tau$. Proof of this:
\begin{eqnarray}
\hat{\tau} & = & L_{\hat{\eta}}\phi=L_{\eta}\phi=L_{\eta}(\chi+g(\tau)[a(\tau)\psi_{\eta_0}+(1-a(\tau))\psi_{\eta}])\nonumber\\
& = & L_{\eta}\chi
+g(\tau)[a(\tau)L_{\eta}\psi_{\eta_0}+(1-a(\tau))L_{\eta}\psi_{\eta}]\nonumber\\
& = & \sigma+(ga)(\tau)L_{\eta}\psi_{\eta_0}\quad\mbox{(with Eq. (4))}\nonumber\\
& = & h_{\eta}(\tau)+(ga)(\tau)L_{\eta}\psi_{\eta_0}=\tau.\nonumber
\end{eqnarray}

2. In case $d(L\chi)=0$ use the result of Part 1 above and continuity as in Part 2 of the proof of Proposition 3. $\Box$

\begin{corollary}
The map $A$ is an almost graph diffeomorphism associated with $X_f$ and $X_0$.
\end{corollary}

{\bf Proof.} Propositions 3 and 4 yield that $A$ is a diffeomorphism onto $C^1$ with inverse $B$. Using $A(X_f)\subset X_0$ and $B(X_0)\subset X_f$ one finds $A(X_f)=X_0$. For $\phi\in X_f\cap X_0$, $\tau=L_{\eta}\phi$, and $\eta=L\phi$, we have  
$$
g(\tau)=g(L_{\eta}\phi)=g(\phi(-d(\eta)))=g(\phi(-d(L\phi)))=\phi'(0)=0.
$$
This yields $A(\phi)=\phi$.
$\Box$

\end{document}